\documentclass[12pt,british]{article}
\usepackage{amsmath,amsfonts,amssymb,amscd}
\usepackage{latexsym,epsfig}
\usepackage{babel}
\usepackage{theorem}
\usepackage[latin1]{inputenc}
\usepackage[T1]{fontenc}
\usepackage{enumerate}
\usepackage{amsbsy}
\usepackage[matrix,arrow]{xy}

\renewcommand{\epsilon}{\ensuremath{\varepsilon}}
\renewcommand{\phi}{\ensuremath{\varphi}}
\renewcommand{\to}{\ensuremath{\longrightarrow}}
\newcommand{\N}{\ensuremath{\mathbb N}}
\newcommand{\Z}{\ensuremath{\mathbb Z}}
\newcommand{\D}{\ensuremath{\mathbb D}}
\newcommand{\T}{\ensuremath{\mathbb{T}^2}}
\newcommand{\FF}{\ensuremath{\mathbb F}}
\newcommand{\F}[1][n]{\ensuremath{\FF_{{#1}}}}

\renewcommand{\ker}[1]{\ensuremath{\operatorname{\text{Ker}}({#1})}}

\makeatletter

\renewcommand{\p@enumii}{}

\def\@enum@{\list{\csname label\@enumctr\endcsname}%
           {\usecounter{\@enumctr}\def\makelabel##1{
\normalfont\ignorespaces\emph{{##1}~}}
\setlength{\labelsep}{3pt}
\setlength{\parsep}{0pt}
\setlength{\itemsep}{5pt}
\setlength{\leftmargin}{0pt}
\setlength{\labelwidth}{0pt}
\setlength{\listparindent}{\parindent}
\setlength{\itemindent}{0pt}
\setlength{\topsep}{3pt plus 1pt minus 1 pt}}}

\def\@map#1#2[#3]{\mbox{$#1 \colon\thinspace #2 \to #3$}}
\def\map#1#2{\@ifnextchar [{\@map{#1}{#2}}{\@map{#1}{#2}[#2]}}

\DeclareRobustCommand*\textsubscript[1]{\@textsubscript{\selectfont#1}}
\def\@textsubscript#1{{\m@th\ensuremath{_{\mbox{\fontsize\sf@size\z@#1}}}}}
\makeatother

\newcommand{\brak}[1]{\ensuremath{\left\{ #1 \right\}}}
\newcommand{\ang}[1]{\ensuremath{\langle #1\rangle}}

\theoremheaderfont{\normalfont\scshape}

\newtheorem{thm}{Theorem}
\newtheorem{lem}[thm]{Lemma}
\newtheorem{prop}[thm]{Proposition}

\theorembodyfont{\rmfamily}

\newcommand{\reth}[1]{Theorem~\protect\ref{th:#1}}
\newcommand{\relem}[1]{Lemma~\protect\ref{lem:#1}}

\newcommand{\resec}[1]{Section~\protect\ref{sec:#1}}

\newcommand{\req}[2][{}]{equation~(\protect\ref{eq:#2}\textsubscript{${#1}$})}

\newcommand{\eop}{%
  \relax
  \ifvmode
    \noindent
  \else
    \unskip
    \hskip0pt plus-1fill\relax 
  \fi
  \vrule width0pt
  \nobreak
  \hfill 
  {\hspace*{\fill}\vrule width3pt height8pt depth0pt}%
}

\newenvironment{proof}{\par\vspace{\partopsep}\noindent\emph{Proof.}}
{\eop\par\vspace{\parsep}}
\newenvironment{prooftext}[1]{\par\vspace{\parsep}\noindent{\emph{Proof of #1.}}}{\eop\par\vspace{\partopsep}}

\newtheorem{rem}[thm]{Remark}
\newtheorem{rems}[thm]{Remarks}

\newcommand{\MG}{\mathcal{M}_g}

\newcommand{\PMG}{\mathcal{M}_{g,p}^{(n)}}
\newcommand{\PPMG}{\mathcal{PM}_{g,p}^{(n)}}
\newcommand{\PPMGa}{\mathcal{PM}_{g,1}^{(n)}}

\begin{document}

\title{Lower central series for surface braid groups}
\author{Paolo Bellingeri\\ 
Laboratoire de Math\'ematiques Jean Leray, UMR~CNRS~6629,\\
Universit\'e de Nantes,\\
2, Rue de la Houssini\`ere,
44072 Nantes, France.\\
bellingeri-p-at-univ-nantes.fr\\
\vspace{5pt}\\
Sylvain Gervais\\
Laboratoire de Math\'ematiques Jean Leray, UMR~CNRS~6629,\\
Universit\'e de Nantes,\\
2, Rue de la Houssini\`ere,
44072 Nantes, France.\\
gervais-s-at-univ-nantes.fr\\
\vspace{5pt}\\
John Guaschi\\
Laboratoire de Math\'ematiques Emile Picard, UMR~CNRS~5580, \\ 
UFR-MIG, Universit\'e Toulouse~III,\\ 
31062~Toulouse Cedex~9, France. \\
guaschi-at-picard.ups-tlse.fr}
\date{\empty}
\maketitle
\begin{abstract}
We determine the lower central series and corresponding residual properties for braid groups and pure braid 
groups of orientable surfaces.
\end{abstract}

\begingroup
\renewcommand{\thefootnote}{}%Removing the footnote symbol.
\endgroup 

\section{Introduction}

Surface braid groups are a natural generalisation of the classical
braid groups (corresponding to the case where $\Sigma$ is a disc) and
of fundamental groups of surfaces (corresponding to the case $n=1$).
They were first defined by Zariski during the 1930's (braid groups on
the sphere had been considered earlier by Hurwitz), were re-discovered
by Fox during the 1960's, and were used subsequently in the study of
mapping class groups.

We recall two definitions of surface braid groups. In
Section~\ref{sec:bound}, we shall give a third equivalent definition
using mapping class groups. 

\paragraph{Surface braid groups via configuration space.} Let $\Sigma$
be a connected, orientable surface. Let $F_n(\Sigma)=\Sigma^n \setminus \Delta$,
where $\Delta$ is the fat diagonal, \emph{i.e.}\ the set of $n$-tuples
$x=(x_1, \dots, x_n)$ for which $x_i=x_j$ for some $i \not= j$. The
fundamental group $\pi_1(F_n(\Sigma))$ is called the \emph{pure braid
group} on $n$ strands of the surface $\Sigma$; it shall be denoted by
$P_n(\Sigma)$. There is a natural action of the symmetric group $S_n$
on $F_n(\Sigma)$ by permutation of coordinates. We denote by 
$\widehat{F_n(\Sigma)}$ the quotient space $F_n(\Sigma)/S_n$. The fundamental group 
$\pi_1(\widehat{F_n(\Sigma)})$ is called the \emph{braid group} on $n$
strands of the surface $\Sigma$; it shall be denoted by $B_n(\Sigma)$.

\paragraph{Surface braid groups as equivalence classes of geometric
braids.} Let $\mathcal{P}=\{p_1, \dots,
p_n\}$ be a set of $n$ distinct points (\emph{punctures}) in the
interior of $\Sigma$. A \emph{geometric braid} on $\Sigma$ based at $\mathcal{P}$
is a collection $(\psi_1, \dots, \psi_n)$ of $n$ disjoint paths
(called \emph{strands}) on $\Sigma \times [0, 1]$ which run
monotonically with $t \in [0, 1]$ and such that
\linebreak[4]$\psi_i(0)=(p_i, 0)$ and $\psi_i(1) \in \mathcal{P}
\times \brak{1}$. Two braids are considered to be equivalent if they
are isotopic. The usual product of paths defines a group structure on
the equivalence classes of braids. This group, which is isomorphic to
$B_n(\Sigma)$, does not depend on the choice of $\mathcal{P}$. A braid
is said to be \emph{pure} if $\psi_i(1)= (p_i, 1)$ for all $i=1,
\dots, n$. The set of pure braids form a group isomorphic to
$P_n(\Sigma)$.

\paragraph{} Given  a group $G$, we define the \emph{lower central series} of $G$
inductively as follows:  set $\Gamma_1(G)=G$, and for $i\geq 2$,
define $\Gamma_i(G)=[G,\Gamma_{i-1}(G)]$. The group $G$ is said to be
\emph{perfect} if $G=\Gamma_2(G)$. From  the lower central series of
$G$ one can define  another  filtration $D_1(G) \supseteq  D_2(G)
\supseteq \cdots$  by setting $D_1(G)=G$, and for $i\geq 2$,
defining $D_i(G)=\{ \, x \in G \, | \, x^n \in \Gamma_{i}(G) $ for some
$n\in \N^* \, \}$. Following Garoufalidis and Levine \cite{GLe}, this filtration is
called  the \emph{rational lower central series} of $G$.\\
Following P.~Hall, for any
group-theoretic property $\mathcal{P}$, a group $G$ is said to be
\emph{residually $\mathcal{P}$} if for any (non-trivial) element $x
\in G$, there exists a group $H$ with the property  $\mathcal{P}$ and
a surjective homomorphism $\map{\phi}{G}[H]$ such that  $\phi(x) \not=1$. It is well known that a group $G$
is residually nilpotent if and only if
$\bigcap_{i \ge 1}\Gamma_i(G)=\{ 1\}$. On the other hand, a group $G$
is residually torsion-free nilpotent if and only if
$\bigcap_{i \ge 1}D_i(G)=\{ 1\}$.

This paper deals with combinatorial properties of surface braid
groups, in particular, their lower central series, and their related
residual properties. In the case of the disc $\mathbb{D}^2$ we have
that  $B_n(\mathbb{D}^2)$ is residually nilpotent if and only if
$n\leq 2$,  and if $n\ge 3$  then $\Gamma_3(B_n(\mathbb{D}^2))=
\Gamma_2(B_n(\mathbb{D}^2))$ (see Proposition~\ref{braidlcs}).
Moreover, Gorin and Lin~\cite{GL} showed that
$\Gamma_2(B_n(\mathbb{D}^2))$  is perfect for $n\ge 5$.

The case of the sphere $\mathbb{S}^2$ and the punctured sphere has been studied by one of the authors and D.~Gon\c{c}alves~\cite{GG2}: 
in particular $B_n(\mathbb{S}^2)$ is residually nilpotent if and only
if $n\leq 2$ and for all $n\ge 3$, $\Gamma_3(B_n(\mathbb{S}^2))= \Gamma_2(B_n(\mathbb{S}^2))$.

Our main results, which concern orientable surfaces of genus at least one, are as follows.

\begin{thm}\label{th:gam3}
Let $\Sigma_g$ be a compact, connected orientable surface without boundary, of genus $g\geq 1$, and let $n\ge
3$. Then:
\begin{enumerate}[(a)]
\item\label{it:gam1g}  $\Gamma_1(B_n(\Sigma_g))/\Gamma_2(B_n(\Sigma_g)) \cong
  \Z^{2g} \oplus \Z_2$. 
\item\label{it:gam2g}  $\Gamma_2(B_n(\Sigma_g))/\Gamma_3(B_n(\Sigma_g)) \cong \Z_{n-1+g}$. 
\item\label{it:gam3g}  $\Gamma_3(B_n(\Sigma_g))=\Gamma_4(B_n(\Sigma_g))$. Moreover
$\Gamma_3(B_n(\Sigma_{g}))$ is perfect for $n\ge 5$.
\item\label{it:gam4g} $B_n(\Sigma_g)$ is not residually nilpotent.
\end{enumerate}
\end{thm}

This implies that braid groups of compact, connected orientable surfaces without boundary may be
distinguished by their lower central series (indeed by the first two lower central quotients).

\begin{thm}\label{th:ggeq2}
Let $g\geq 1$, $m\geq 1$ and $n\geq 3$. Let $\Sigma_{g,m}$ be a
compact, connected orientable surface of genus $g$ with $m$ boundary
components. Then:
\begin{enumerate}[(a)]
\item\label{it:g1sigma} $\Gamma_1(B_n(\Sigma_{g,m}))/ \Gamma_2(B_n(\Sigma_{g,m}))=\Z^{2g+m-1} \oplus \Z_2$.
\item\label{it:g2sigma} $\Gamma_2(B_n(\Sigma_{g,m}))/ \Gamma_3(B_n(\Sigma_{g,m}))=\Z$.
\item\label{it:g3sigma} $\Gamma_3(B_n(\Sigma_{g,m}))= \Gamma_4(B_n(\Sigma_{g,m}))$. Moreover
$\Gamma_3(B_n(\Sigma_{g,m}))$ is perfect for $n\ge 5$.
\item\label{it:resid} $B_n(\Sigma_{g,m})$ is not residually nilpotent.
\end{enumerate}
\end{thm}

Braid groups on $2$ strands represent a very difficult and interesting
case. In the case of the torus, we are able to prove that its $2$-strand braid group is residually nilpotent.
Further, using ideas from~\cite{GG2} and results of~\cite{Ga}, we may show that apart from the first term,
the lower central series of $B_2(\T)$ and $\Z_2 \ast \Z_2 \ast \Z_2$ coincide, and we may also determine all
of their successive lower central quotients. More precisely:

\begin{thm}\label{th:renil}\mbox{}
\begin{enumerate}%[(a)]
\item $B_2(\T)$ is residually nilpotent.
\item\label{it:igeq2} For all $i\geq 2$:
\begin{enumerate}%[(i)]
\item\label{it:gammai} $\Gamma_i(B_2(\T))\cong \Gamma_i(\Z_2 \ast \Z_2 \ast \Z_2)$.
\item $\Gamma_i(B_2(\T))/\Gamma_{i+1}(B_2(\T))$ is isomorphic to the direct sum of $R_i$ copies of $\Z_2$,
where:
\begin{equation*}
R_i=\sum_{j=1}^{i-2} \left( \sum_{
\substack
{k \mid i-j\\
k>1}}
\; \mu \left( \frac{i-j}{k}\right) \frac{k\alpha_k}{i-j}
\right) \quad \text{and} \quad k\alpha_k=2^k+2(-1)^k.
\end{equation*}
\end{enumerate}
\item $B_2(\T)$ is not residually torsion-free nilpotent.
\end{enumerate}
\end{thm}

Finally, as we shall see in
Proposition~\ref{biord}, $B_2(\T)$ is not bi-orderable (see \resec{at}
for a definition).

In Section~\ref{puresection}, we recall the relations between mapping
class groups and surface braid groups, and prove that  pure braid
groups of the torus and of surfaces with boundary components are
residually torsion-free nilpotent. This is achieved by showing that
they may be realised as  subgroups of the Torelli group of a surface
of higher genus (Lemma~\ref{lem:puretorelli}), which is known to be
residually torsion-free nilpotent (see for instance~\cite{H}). We note that the embedding proposed in Lemma~\ref{lem:puretorelli} 
does not hold when the surface is without boundary (see Remark~\ref{rem:closed}).

\vspace{5pt}

\noindent \textbf{Acknowledgements} The authors are grateful to
I.~Marin for pointing out the relevance of the notion of residual
torsion-free nilpotence, to V.~Bardakov for useful suggestions,
and to L.~Paris for the reference~\cite{H}.

%%%%%%%%%%%%%%%%%%%%%%%%%%%%%%%%%%%%
%%%%%%%%%%%%%%%%%%%%%%%%%%%%%%%%%%%%
%%%%%%%%%%%%%%%%%%%%%%%%%%%%%%%%%%
%%%%%%%%%%%%%%%%%%%%%%%%%%%%%%%%%%%%
%%%%%%%%%%%%%%%%%%%%%%%%%%%%%%%%%%%%
%%%%%%%%%%%%%%%%%%%%%%%%%%%%%%%%%%
%%%%%%%%%%%%%%%%%%%%%%%%%%%%%%%%%%%%
%%%%%%%%%%%%%%%%%%%%%%%%%%%%%%%%%%%%
%%%%%%%%%%%%%%%%%%%%%%%%%%%%%%%%%%
%%%%%%%%%%%%%%%%%%%%%%%%%%%%%%%%%%%%
%%%%%%%%%%%%%%%%%%%%%%%%%%%%%%%%%%%%
%%%%%%%%%%%%%%%%%%%%%%%%%%%%%%%%%%
%%%%%%%%%%%%%%%%%%%%%%%%%%%%%%%%%%%%
%%%%%%%%%%%%%%%%%%%%%%%%%%%%%%%%%%%%

\section{Lower central series for Artin-Tits groups}\label{sec:at}

Let us start by recalling some standard results on combinatorial
properties of braid groups.  The following result is well known
(see~\cite{GL} for instance).

\begin{prop} \label{braidlcs}
Let $B_n$ be the Artin braid group on $n\ge 3$ strands.\\
Then $\Gamma_1(B_n)/\Gamma_2(B_n) \cong \Z$ and $\Gamma_2(B_n)=\Gamma_3(B_n)$.
\end{prop}

\begin{proof}
Let us give an easy proof of the second statement (we use an
argument of~\cite{GG2}). Let $\{ \sigma_1, \dots, \sigma_{n-1} \}$ be
the usual set of generators of $B_n$; the classical relations of
$B_n$, referred to hereafter as \emph{braid relations}, are as
follows:
\begin{gather}
\text{$\sigma_i\sigma_j=\sigma_j\sigma_i$, for all $1\leq i,j\leq n-1$ and $\lvert i-j\rvert \geq
2$,}\label{eq:artinbn1}\\
\text{$\sigma_i\sigma_{i+1}\sigma_i=\sigma_{i+1} \sigma_i\sigma_{i+1}$ for all $1\leq i\leq
n-2$.}\label{eq:artinbn2}
\end{gather}
From this, we see that $B_n/\Gamma_2(B_n)$ is isomorphic to $\Z$.

Consider the following short exact sequence:
\begin{equation*}
\displaystyle{1 \to \frac{\Gamma_2(B_n)}{\Gamma_3(B_n)} \to
  \frac{B_n}{\Gamma_3(B_n)} \stackrel{p}{\to} \frac{B_n}{\Gamma_2(B_n)}\to 1},
\end{equation*} 
Since all of the $\sigma_i\in B_n/\Gamma_3(B_n)$ project to the same
element of $B_n/\Gamma_2(B_n)$, for each $1\leq i\leq n-1$, there
exists $t_i\in \Gamma_2(B_n)/\Gamma_3(B_n)$ (with $t_1=1$) such that
$\sigma_i=t_i\sigma_1$. Projecting the braid
relation~(\ref{eq:artinbn2})  into $B_n/\Gamma_3(B_n)$, we see that
$t_i \sigma_1 t_{i+1}\sigma_1 t_i \sigma_1= t_{i+1} \sigma_1
t_i\sigma_1 t_{i+1} \sigma_1$. But the $t_i$ are central in $B_n /
\Gamma_3(B_n)$, so $t_i=t_{i+1}$, and since $t_1=1$, we obtain
$\sigma_1=\ldots =\sigma_{n-1}$. So the surjective homomorphism $p$ is
in fact an isomorphism. 
\end{proof}

We recall that classical braid groups are also called Artin-Tits
groups of type $\mathcal{A}$.  More precisely, let  $(W,S)$ be 
a Coxeter system and let us denote 
by $m_{s, \, t}$ the order of the element $s t$ in $W$ (for $s, t \in S$).
Let $B_W$ be the group defined by the following group presentation:
$$
B_W=\langle S \, | \underbrace{s t\cdots}_{m_{s, \, t}}= \underbrace{t s\cdots}_{m_{s, \, t}}\; \mbox{ for any $s\not=t \in S$
with $m_{s, \, t} < +\infty$} \, \rangle
$$
The group $B_W$ is the Artin-Tits group associated to $W$. The group
$B_W$ is said to be of spherical type if  $W$ is finite. The kernel of
the canonical projection of $B_W$ onto $W$ is called the pure Artin-Tits group associated to $W$.

\begin{prop} Let $B_W$ be an Artin-Tits group of spherical type, with $W$ different from the dihedral group 
$I_{2m}$. Then:
\begin{enumerate}[i)]
\item $\Gamma_1(B_W)/\Gamma_2(B_W)$ is isomorphic to $\Z$ or $\Z^2$.
\item  $\Gamma_2(B_W)=\Gamma_3(B_W)$.
\end{enumerate}
\end{prop}
\begin{proof}
Considering  $B_W$  equipped with the previous group presentation,  it is easy to calculate the Abelianisation.
Using the same argument as in Proposition~\ref{braidlcs}, one deduces that
if $s$ and $t$ are two generators of an Artin-Tits group of spherical type $B_W$ then they are identified in 
$B_W/\Gamma_3(B_W)$ if $m_{s, \, t}$ is odd.
This argument allows us to prove that $B_W/\Gamma_3(B_W)$ is isomorphic to 
$B_W/\Gamma_2(B_W)$ for almost all Artin-Tits groups of spherical type, the only exception being that of the one-relator group
$B_{I_{2m}}=\langle a,b \, | \, (ab)^m=(ba)^m \rangle$ (for $m>1$),
since the defining relation is of even length.
\end{proof} 

\begin{rems}
The Artin-Tits group $B_{I_{2m}}=\langle a,b \, | \, (ab)^m=(ba)^m \rangle$
 is residually nilpotent 
if and only if $m$ is a power of a prime number. Indeed, by taking $c=ba$, it is readily seen that the group $B_W$ is isomorphic to
the Baumslag-Solitar group of type $(m,m)$, $BS_m=\langle a,c \, | \, [a, \, c^m]=1 \rangle$,
which is known to  be residually nilpotent if $m$ is a power of a prime number. Conversely, let $G$ be a one-relator group with non-trivial centre. According to \cite{McC}, $G$ is residually nilpotent if and only if one of the following holds:
\begin{enumerate}[i)]
\item 
$G$ is Abelian.
\item $G$ is
isomorphic to a Baumslag-Solitar group of type $(r,r)$, with  $r$ a power of a prime number.
\item
$G$ is isomorphic to $G_{p,q}=\langle m, n \, | \, m^p=n^q \rangle$ 
where $p$ and $q$ are powers of the same prime number.
\end{enumerate}
Suppose that $B_{I_{2m}}$ is residually nilpotent. The group
$B_{I_{2m}}$ is not Abelian for $m>1$, and it cannot be isomorphic to
a group $G_{p,q}$ since they  have different Abelianisations.
Therefore  $B_{I_{2m}}$ is isomorphic to a Baumslag-Solitar group of type
$(r,r)$ with $r$ a power of a prime number. But as we remarked,
$B_{I_{2m}}$ is also isomorphic to the Baumslag-Solitar group of type
$(m,m)$, in which case $m=r$, by a result of Moldavanskii on
isomorphisms of Baumslag-Solitar groups~(\cite{Mol}). It thus follows that $B_{I_{2m}}$ is residually nilpotent if and only if $m$ is a power of a prime number.
\end{rems}

On the other hand, it is well known that pure braid groups are residually torsion-free nilpotent~\cite{FR}. Using the faithfulness of the
Krammer-Digne
representation, Marin has shown recently that the pure Artin-Tits
groups of spherical type  are  residually torsion-free nilpotent~\cite{M}.

The fact that a group is residually torsion-free nilpotent has several
important consequences, notably that the group is bi-orderable~\cite{MR}. We recall that a group $G$ is said
to be
\emph{bi-orderable} if there exists a strict total ordering $<$ on its
elements which is invariant under left and right multiplication, in
other words, $g<h$ implies that $gk<hk$ and  $kg<kh$ for all $g,h,k\in
G$. We state one interesting property of bi-orderable groups. A group
$G$ is said to have \emph{generalised torsion} if there exist $g,
h_1, \dots, h_k$, $(g\not=1)$ such that: 
\begin{equation*}
(h_1 g h_1^{-1}) (h_2 g h_2^{-1}) \cdots (h_k g h_k^{-1})=1 \,.
\end{equation*}

\begin{prop}[\cite{KK}]\label{test}
A bi-orderable group has no generalised torsion. 
\end{prop}

The braid group $B_n$ is not bi-orderable for $n \ge 3$ since it has
generalised torsion (see~\cite{N} or~\cite{ba}).  As we shall see in
Section~\ref{sec:b2t2}, $B_2(\mathbb{T}^2)$ \emph{is} residually
nilpotent, but is not bi-orderable.

%%%%%%%%%%%%%%%%%%%%%%%%%%%%%%%%%%%%%%%%%%%%%%%%%%%%%%%%%%%%%%%%%%%%%%%
%%%%%%%%%%%%%%%%%%%%%%%%%%%%%%%%%%%%%%%%%%%%%%%%%%%%%%%%%%%%%%%%%%%%%%%
%%%%%%%%%%%%%%%%%%%%%%%%%%%%%%%%%%%%%%%%%%%%%%%%%%%%%%%%%%%%%%%%%%%%%%%
%%%%%%%%%%%%%%%%%%%%%%%%%%%%%%%%%%%%%%%%%%%%%%%%%%%%%%%%%%%%%%%%%%%%%%%
%%%%%%%%%%%%%%%%%%%%%%%%%%%%%%%%%%%%%%%%%%%%%%%%%%%%%%%%%%%%%%%%%%%%%%%
%%%%%%%%%%%%%%%%%%%%%%%%%%%%%%%%%%%%%%%%%%%%%%%%%%%%%%%%%%%%%%%%%%%%%%%
%%%%%%%%%%%%%%%%%%%%%%%%%%%%%%%%%%%%%%%%%%%%%%%%%%%%%%%%%%%%%%%%%%%%%%%
%%%%%%%%%%%%%%%%%%%%%%%%%%%%%%%%%%%%%%%%%%%%%%%%%%%%%%%%%%%%%%%%%%%%%%%
%%%%%%%%%%%%%%%%%%%%%%%%%%%%%%%%%%%%%%%%%%%%%%%%%%%%%%%%%%%%%%%%%%%%%%%

\section{Lower central series for surface braid groups on at least $3$ strands}

\subsection{Surfaces without boundary}

This section is devoted to proving \reth{gam3}. Let $\Sigma_g$ be a compact, connected orientable surface
without boundary, of genus $g>0$. We start by giving a presentation of $B_n(\Sigma_g)$.

\begin{thm}(\cite{B})\label{th:presbng}
Let $n\in\N$. Then $B_n(\Sigma_g)$ admits the following group presentation:
\begin{enumerate}
\item[\textbf{Generators:}] $a_1, b_1, \ldots, a_g, b_g, \sigma_1, \ldots, \sigma_{n-1}$.
\item[\textbf{Relations:}] 
\begin{gather}
\text{$\sigma_i\sigma_j=\sigma_j\sigma_i$ if $\lvert i-j \rvert \geq
2$}\label{eq:artin1}\\
\text{$\sigma_i\sigma_{i+1}\sigma_i= \sigma_{i+1}\sigma_i \sigma_{i+1}$
for all $1\leq i\leq n-2$}\label{eq:artin2}\\
\text{$c_i\sigma_j= \sigma_j c_i$ for all $j\geq 2$,   $c_i=a_i$ or
  $b_i$ and $i=1, \ldots, g$}\label{eq:asjg}\\
\text{$c_i \sigma_1 c_i \sigma_1= \sigma_1 c_i \sigma_1 c_i$   for
  $c_i=a_i$ or $b_i$ and $i=1, \ldots, g$}\label{eq:bsbg}\\
\text{$a_i \sigma_1 b_i = \sigma_1 b_i \sigma_1 a_i \sigma_1$ for  $i=1, \ldots, g$}\label{eq:abbag}\\
\text{$c_i \sigma_1^{-1}  c_j \sigma_1=\sigma_1^{-1} c_j \sigma_1 c_i$
  for $c_i=a_i$ or $b_i$, $c_j=a_j$ or $b_j$ and $1\le j<i\le g$}\label{eq:cddcg}\\
\text{$\prod_{i=1}^g [a_i^{-1},b_i]= \sigma_1\cdots \sigma_{n-2} \sigma_{n-1}^2 \sigma_{n-2} \cdots \sigma_1$}\label{eq:abcommg}.
\end{gather}
\end{enumerate}
\end{thm}
\begin{proof}
Let $\widetilde{B_n(\Sigma_g)}$ be the group defined by the above
presentation, and let $B_n(\Sigma_g)$ be the group given by the
presentation of Theorem~1.2 of~\cite{B}.  Consider the homomorphism
$\map{\phi}{B_n(\Sigma_g)}[\widetilde{B_n(\Sigma_g)}]$ 
defined on the
generators of $B_n(\Sigma_g)$ by $\phi(\sigma_j)=\sigma_j$ (for $j=1,
\ldots, n-1$), $\phi(a_i)=a_i^{-1}$ and $\phi(b_i)=b_i^{-1}$ (for
$i=1, \ldots, g$).  It is an easy exercise to check that $\phi$ is an
isomorphism. 
\end{proof}

\begin{prooftext}{\reth{gam3}}
\begin{enumerate}[(a)]
\item Consider the group $\Z^{2g} \oplus \Z_2$ defined by the
presentation  $\langle c_1, \dots, c_{2g},  \sigma \, | \, \sigma^2=
[c_i, c_j]= [c_i, \sigma]=1, \,$ for $\, 1 \le i,j\le 2g \,\rangle$
and  $B_n(\Sigma_g)$ with the group presentation given by
\reth{presbng}. It is easy to check  that the homomorphism
\begin{equation*}
\map{\phi}{\Gamma_1(B_n(\Sigma_g))/\Gamma_2(B_n(\Sigma_g))}[\Z^{2g}
\oplus \Z_2]
\end{equation*}
which sends $a_k$ to $c_{2k-1}$, $b_k$ to $c_{2k}$ and
every $\sigma_j$ to $\sigma$ is indeed an isomorphism.

\item Let us start by determining a group presentation for
$B_n(\Sigma_g)/\Gamma_3(B_n(\Sigma_g))$. Let $q$ be the canonical 
projection of $B_n(\Sigma_g)$ onto
$B_n(\Sigma_g)/\Gamma_3(B_n(\Sigma_g))$. As in the proof of
Proposition~\ref{braidlcs}, the braid relations~(\ref{eq:artin2})
imply that $q(\sigma_1)= \dots= q(\sigma_{n-1})$; we denote this
element by $\sigma$.  This implies that the projected
relations~(\ref{eq:artin1}) are trivial. For $i=1, \dots, g$, let us
also denote $q(a_i)$ by $a_i$ and $q(b_i)$ by $b_i$. Since $n\ge 3$,
we see from relations~(\ref{eq:asjg}) that $\sigma$ is central in 
$B_n(\Sigma_g)/\Gamma_3(B_n(\Sigma_g))$ and hence the projected
relations~(\ref{eq:bsbg}) become trivial. From
relations~(\ref{eq:cddcg}), for all $1\leq i,j\leq g$, $i \not= j$,
one may infer that $[a_i, b_j]=[a_i, a_j]=[b_i, b_j] =[b_i, a_j]=1$ in
$B_n(\Sigma_g)/\Gamma_3(B_n(\Sigma_g))$.   Relations~(\ref{eq:abbag})
and~(\ref{eq:abcommg}) imply that $[b_i,a_i]=\sigma^{-2}$ for all $=1,
\dots, g$, and $\prod_{i=1}^g [a_i^{-1},b_i]= \sigma^{2(n-1)}$
respectively. Conjugating the latter equation by $a_1\cdots a_g$
yields $\prod_{i=1}^g [b_i,a_i]= \sigma^{2(n-1)}$ in
$B_n(\Sigma_g)/\Gamma_3(B_n(\Sigma_g))$ (recall that $a_i$ commutes
with $a_j$ and $b_j$ if $i\not=j$), and hence
$\sigma^{-2g}=\prod_{i=1}^g [b_i, a_i]= \sigma^{2(n-1)}$. Therefore,
$\sigma^{2(g+n-1)}=1$ in $B_n(\Sigma_g)/\Gamma_3(B_n(\Sigma_g))$.
Summing up, we have obtained the following information:
\begin{equation}\label{eq:presbngam4}
\left.\begin{aligned}
&\text{$B_n(\Sigma_g)/\Gamma_3(B_n(\Sigma_g))$ is generated by
$a_1,b_1, \ldots, a_g, b_g$ and $\sigma$}\\
&\text{$a_1,b_1, \ldots, a_g, b_g$ and $\sigma$ commute pairwise except for the pairs $(a_i, b_i)_{i=1,
\ldots, g}$}\\
&[a_1,  b_1]=\cdots=   [a_g,  b_g] =\sigma^2; \qquad \sigma^{2(n+g-1)}=1.
\end{aligned} \right\}
\end{equation}

The remaining relations of $B_n(\Sigma_g)/\Gamma_3(B_n(\Sigma_g))$ are
those of the form $[[x,y],z]=1$  for all $x,y,z\in
B_n(\Sigma_g)/\Gamma_3(B_n(\Sigma_g))$. We claim that such relations
are implied by those of~(\ref{eq:presbngam4}).  To see this, recall
that $\Gamma_2(B_n(\Sigma_g)/\Gamma_3(B_n(\Sigma_g)))$ is the normal
subgroup of  $B_n(\Sigma_g)/\Gamma_3(B_n(\Sigma_g))$ generated by the
finite set of commutators $[a_i, a_j]$,  $[b_i, b_j]$,  $[a_i, b_j]$,
$[a_i, b_i]$, $[a_i, \sigma]$ and $[b_i, \sigma]$, for $1 \le i \not=
j \le g$. But the relations of~(\ref{eq:presbngam4}) imply that these
commutators are all trivial, with the exception of $[a_i, b_i]$ for $1
\le i \le g$, which is equal to $\sigma^2$.  Since $\sigma$ is central
in $B_n(\Sigma_g)/\Gamma_3(B_n(\Sigma_g))$, we conclude that
$\Gamma_2(B_n(\Sigma_g)/\Gamma_3(B_n(\Sigma_g)))= \ang{\sigma^2}$, and
that $[[x,y],z]=1$ as claimed. Hence~(\ref{eq:presbngam4}) is a group
presentation for $B_n(\Sigma_g)/\Gamma_3(B_n(\Sigma_g))$.

Now consider the following exact sequence:
\begin{equation*}
\displaystyle{1 \to \frac{\Gamma_2(B_n(\Sigma_g))}{\Gamma_3(B_n(\Sigma_g))} \to
  \frac{B_n(\Sigma_g)}{\Gamma_3(B_n(\Sigma_g))} \stackrel{p}{\to}
\frac{B_n(\Sigma_g)}{\Gamma_2(B_n(\Sigma_g)}\to 1}.
\end{equation*}
From the presentation of $B_n(\Sigma_g)/\Gamma_3(B_n(\Sigma_g))$
given by~(\ref{eq:presbngam4}), one sees that every element $x$ of 
$B_n(\Sigma_g)/\Gamma_3(B_n(\Sigma_g))$ may be written in the form
$a_1^{j_1} b_1^{k_1} \cdots a_g^{j_g} b_g^{k_g} \sigma^p$. Since
$B_n(\Sigma_g)/\Gamma_2(B_n(\Sigma_g))$  is isomorphic to $\Z^{2g} \oplus
\Z_2$, the factors being generated respectively by  $p(a_1),p(b_1),
\ldots, p(a_g), p(b_g)$ and $p(\sigma)$, if $x\in \ker p$ then
$j_1=k_1=\ldots=j_g=k_g=0$ and $p$ is even, so $\ker p\subseteq
\ang{\sigma^2}$. The converse is clearly true and so $\ker
p=\ang{\sigma^2}$. 

Let $d$ denote the order of $\sigma^2$ in $B_n(\Sigma_g)/\Gamma_3(B_n(\Sigma_g))$.
From~(\ref{eq:presbngam4}), we have that $d$ divides $n+g-1$.  To complete the proof of part~(\ref{it:gam2g})
of \reth{gam3}, it
suffices to show that $n+g-1$ divides $d$.

Let $G$ be the group generated by elements $a_1, b_1, \ldots, a_g,
b_g$ and $\sigma$, whose relations are $\sigma^{2(n+g-1)}=1$, and the
generators commute pairwise except for the pairs $(a_i, b_i)$ for
$i=1, \ldots, g$. Then $G=\left( \oplus_{i=1}^{g}\F[2](a_i,b_i)
\right) \oplus \Z_{2(n+g-1)}$, where $\F[2](a_i,b_i)$ denotes the free
group of rank~$2$ generated by $a_i$ and $b_i$. Let $N$ be the
subgroup of $G$ normally generated by the elements
$[a_1,b_1]\sigma^{-2}, \ldots, [a_g, b_g]\sigma^{-2}$, and let $\rho$
denote the canonical projection $G\to G/N$.  Then $G/N\cong
B_n(\Sigma_g)/\Gamma_3(B_n(\Sigma_g))$ by the group presentation given
in (\ref{eq:presbngam4}). The cosets modulo~$N$ of the elements $a_1,
b_1, \ldots, a_g, b_g$ and $\sigma$ of $G$ may be identified
respectively with the elements $a_1, b_1, \ldots, a_g, b_g$ and
$\sigma$ of $B_n(\Sigma_g)/\Gamma_3(B_n(\Sigma_g))$.  Further, by
applying the relations of $G$, any element of  $N$ may be written in
the form  $\prod_{i=1}^g \left( \prod_{k=1}^{m_i}\; u_{i_k}
[a_i,b_i]^{\epsilon_{i_k}}u_{i_k}^{-1} \right) \sigma^{-2
\left(\sum_{i=1}^g \left( \sum_{k=1}^{m_i} \; \epsilon_{i_k} \right)
\right) }$, where $m_i \in \N$ for all $i=1, \ldots, g$, and for all
$k=1, \dots, m_i$, $u_{i_k}\in \F[2](a_i,b_i)$ and $\epsilon_{i_k} \in
\brak{1,-1}$. Since $\sigma^{2d}=1$ in
$B_n(\Sigma_g)/\Gamma_3(B_n(\Sigma_g))$, and so in $G/N$, considered
as an element of $G$, it follows that $\sigma^{2d}$ belongs to
$\ker{\rho}$. Hence for all $i=1, \ldots, g$, there exists  $m_i \in
\N$, and for $1\leq k \leq m_i$, there exist $u_{i_k}\in
\F[2](a_i,b_i)$ and $\epsilon_{i_k} \in \brak{1,-1}$ such that  
$\sigma^{2d}=\prod_{i=1}^g \left( \prod_{k=1}^{m_i}\;
u_{i_k}[a_i,b_i]^{\epsilon_{i_k}}u_{i_k}^{-1} \right) \sigma^{-2
\left(\sum_{i=1}^g \left( \sum_{k=1}^{m_i} \; \epsilon_{i_k} \right)
\right)}$. Thus:
\begin{equation}\label{eq:sigabg}
\sigma^{2 \left( d+\sum_{i=1}^g \left( \sum_{k=1}^{m_i} \; \epsilon_{i_k}
\right) \right) }=\prod_{i=1}^g
\left( \prod_{k=1}^{m_i}\;
u_{i_k}([a_i,b_i])^{\epsilon_{i_k}}u_{i_k}^{-1} \right).
\end{equation}
From the structure of $G$, it follows that both the right- and
left-hand sides are equal to~$1$. Moreover,
$\Gamma_2(G)=\oplus_{i=1}^g \Gamma_2(\F[2](a_i,b_i))$. Let  $1\leq
i\leq g$. Projecting the right-hand side of \req{sigabg}, which
belongs to $\Gamma_2(G)$, into $\Gamma_2(\F[2](a_i,b_i))$, and then
into $\Gamma_2(\F[2](a_i,b_i))/ \Gamma_3(\F[2](a_i,b_i))$, we observe
that  $[a_i,b_i]^{\sum_{k=1}^{m_i} \; \epsilon_{i_k}}=1$.  But this
quotient is an infinite cyclic group~\cite{MKS}, hence
$\sum_{k=1}^{m_i} \; \epsilon_{i_k}=0$ for $i=1, \ldots, g$ and
therefore $\sum_{i=1}^g \left( \sum_{k=1}^{m_i} \;
\epsilon_{i_k}\right)=0$.  Thus the left-hand side of \req{sigabg}
reduces to $\sigma^{2d}=1$ in $G$, and so $n+g-1$ divides $d$. It
follows that $\sigma$ is of order~$2(n+g-1)$ in
$B_n(\Sigma_g)/\Gamma_3(B_n(\Sigma_g))$ as claimed.

\item Let $H$ denote the normal closure in $B_n(\Sigma_g)$ of the
element $\sigma_1\sigma_2^{-1}$. Using the Artin braid relations, one
may check that in $B_n(\Sigma_g)/H$, the $\sigma_i$ are all identified
to a single element, $\sigma$,  say, and then that \req{presbngam4}
defines a group presentation for $B_n(\Sigma_g)/H$. Thus
$B_n(\Sigma_g)/H\cong B_n(\Sigma_g)/\Gamma_3(B_n(\Sigma_g))$ via  an
isomorphism $\iota$. Now $B_n(\Sigma_g)$ contains a copy of the usual
Artin braid group $B_n$ which is generated by the $\sigma_i$. From the
Artin braid relations, it follows that $\Gamma_2(B_n)$ is the normal
closure in $B_n$ of the elements $\sigma_i \sigma_{i+1}^{-1}$, $1\leq
i\leq n-2$. Moreover, since  $\sigma_{i+1} \sigma_{i+2}^{-1}=
\sigma_i^{-1}\sigma_{i+1}^{-1}\sigma_{i+2}^{-1}\cdot \sigma_i
\sigma_{i+1}^{-1} \cdot \sigma_{i+2} \sigma_{i+1} \sigma_i$ for all
$1\leq i\leq n-3$, we see that $\Gamma_2(B_n)$ is the normal closure 
in $B_n$ of just $\sigma_1 \sigma_2^{-1}$, and thus
$\mathcal{N}_{B_n(\Sigma_g)}(\Gamma_2(B_n))=H$ (if $X$ is a subset of
a group $G$, then we denote its normal closure in $G$ by
$\mathcal{N}_G(X)$).

Since $\Gamma_3(B_n)=\Gamma_2(B_n)$ (by Proposition~\ref{braidlcs}),
we have $\Gamma_4(B_n(\Sigma_g)) \supseteq \Gamma_4(B_n)=
\Gamma_2(B_n)$. Taking normal closures in $B_n(\Sigma_g)$, we deduce
that  $H$ is a normal subgroup of $\Gamma_4(B_n(\Sigma_g))$, and hence
we obtain the following commutative diagram:
\begin{equation*}
\xymatrix{%
1 \ar[r] & \Gamma_4(B_n(\Sigma_g))/H \ar[r] &
B_n(\Sigma_g)/H \ar[drr]_{\cong}^{\iota} \ar@{>>}[rr] & & B_n(\Sigma_g)/\Gamma_4(B_n(\Sigma_g)) \ar@{>>}[d]
\ar[r] & 1\\
& & & & B_n(\Sigma_g)/\Gamma_3(B_n(\Sigma_g)).  &
}
\end{equation*}

Since $\iota$ is an isomorphism, so is the vertical arrow, and hence
its kernel
\linebreak[4]$\Gamma_3(B_n(\Sigma_g))/\Gamma_4(B_n(\Sigma_g))$ is
trivial. This proves the first part of~(\ref{it:gam3g}). To prove the
second part, we have just seen that the normal closure $H$ of
$\Gamma_2(B_n)$ in $B_n(\Sigma_g)$, is isomorphic to
$\Gamma_3(B_n(\Sigma_g))$  (they coincide in fact). Since
$\Gamma_2(B_n)$ is perfect for all $n \ge 5$~\cite{GL}, so are $H$ and
$\Gamma_3(B_n(\Sigma_g))$.

\item We first remark that $\Gamma_3(B_n(\Sigma_g))\neq \brak{1}$. For if $\Gamma_3(B_n(\Sigma_g))$ were
trivial, by~(\ref{it:gam3g}), we would have $\Gamma_2(B_n(\Sigma_g))\cong \Z_{n-1+g}$. But by~\cite{VB},
$B_n(\Sigma_g)$ is torsion free, which yields a contradiction. From this it follows that $\bigcap_{i\in\N}
\Gamma_i(B_n(\Sigma_g))\neq\brak{1}$. This completes the proof of the theorem.
\end{enumerate}
\end{prooftext}

\begin{rem}
Given a group $G$, the property that the $i$th term $\Gamma_i(G)$ is
perfect implies that  $\Gamma_i(G)=\Gamma_{i+1}(G)$.
\end{rem}

%%%%%%%%%%%%%%%%%%%%%%%%%%%%
%%%%%%%%%%%%%%%%%%%%%%%%%%%%
%%%%%%%%%%%%%%%%%%%%%%%%%%%%
%%%%%%%%%%%%%%%%%%%%%%%%%%%%
%%%%%%%%%%%%%%%%%%%%%%%%%%%%
%%%%%%%%%%%%%%%%%%%%%%%%%%%%
%%%%%%%%%%%%%%%%%%%%%%%%%%%%
%%%%%%%%%%%%%%%%%%%%%%%%%%%%
%%%%%%%%%%%%%%%%%%%%%%%%%%%%

\subsection{Surfaces with non-empty boundary}

In this section, we study the case of orientable surfaces with
boundary, and prove \reth{ggeq2}. We identify $\Sigma_{g,0}$ with
$\Sigma_g$. As in \reth{presbng}, from Theorem~1.1 of~\cite{B}, one
obtains the following presentation of $B_n(\Sigma_{g,m})$.

\begin{thm}\label{th:presbnpg}
Let $n\in\N$. Then $B_n(\Sigma_{g,m})$ admits the following group presentation:
\begin{enumerate}
\item[\textbf{Generators:}] $a_1, b_1, \ldots, a_g, b_g, z_1, \dots
  z_{m-1}, \sigma_1, \ldots, \sigma_{n-1}$.
\item[\textbf{Relations:}] 
\begin{gather*}
 \text{Relations~(\ref{eq:artin1})~--~(\ref{eq:cddcg})
 of \reth{presbng}}
\end{gather*}
\vspace{-1 cm}
\begin{gather}
\text{$z_i\sigma_j= \sigma_j z_i$ for all $j\geq 2$ and $i=1, \ldots, m-1$}\label{eq:asjgp}\\
\text{$z_i \sigma_1 z_i \sigma_1= \sigma_1 z_i \sigma_1 z_i$ for  $i=1, \ldots, m-1$}\label{eq:bsbpg}\\
\text{$z_i \sigma_1^{-1}   z_j \sigma_1= \sigma_1^{-1} z_j \sigma_1
  z_i$ for  $1 \le  j < i \le m-1$}\label{eq:asapg}\\
\text{$c_i \sigma_1^{-1}  z_j \sigma_1=\sigma_1^{-1} z_j \sigma_1 c_i$
  for $c_i=a_i$ or $b_i$,  $i=1, \ldots, g$ and  $j=1, \ldots, m-1$.}\label{eq:cddpcg} 
\end{gather}
\end{enumerate}
\end{thm}

\begin{prooftext}{\reth{ggeq2}}
Statement~(\ref{it:g1sigma}) may be proved in the same way as~(\ref{it:gam1g}) of \reth{gam3}.

We now prove part~(\ref{it:g2sigma}). As in the proof of
part~(\ref{it:gam2g}) of \reth{gam3}, one may check that
$\Gamma_2(B_n(\Sigma_{g,m}))/ \Gamma_3(B_n(\Sigma_{g,m}))=\langle
\sigma^2\rangle$, where for all $i=1,\ldots,n-1$, $\sigma$ is the
projection of $\sigma_i$ in $B_n(\Sigma_{g,m})/
\Gamma_3(B_n(\Sigma_{g,m})$. It thus suffices to show that $\sigma^2$
is of infinite order.

Instead of repeating the arguments used in \reth{gam3},  we propose a
different proof, based on geometric relations between surface braid
groups. Suppose that $\sigma^{2d}=1$ for some $d \in \N$. This is
equivalent to saying that $\sigma_i^{2d}$ belongs to
$\Gamma_3(B_n(\Sigma_{g,m}))$ for all $i=1, \dots, n-1$.

Let $1\leq i\leq m$. To each boundary component $\partial_i$ of
$\Sigma_{g,m}$ let us associate a surface $\Sigma_{g_i, 1}$ of
positive genus $g_i$. We choose the $g_i$ so that $h=g+\sum_{i=1}^m
g_i > d-(n-1)$. Let $\Sigma_h$ denote the compact, orientable surface
without boundary and of genus~$h$ obtained by glueing $\partial
\Sigma_{g_i, 1}$ to $\partial_i$ for all $i=1,\ldots, m$. The
embedding of $\Sigma_{g,m}$ into $\Sigma_h$ induces a natural
homomorphism $\phi$ between $B_n(\Sigma_{g,m})$ and $B_n(\Sigma_h)$,
sending geometric generators of $B_n(\Sigma_{g,m})$ to the
corresponding elements of $B_n(\Sigma_h)$. In particular,
$\phi(\sigma_i)=\sigma_i$ for all $i=1, \dots,n-1$.

Since $\sigma_i^{2d}$ belongs to $\Gamma_3(B_n(\Sigma_{g,m}))$, it
follows that $\phi(\sigma_i^{2d})= \sigma_i^{2d}$ belongs to
$\Gamma_3(B_n(\Sigma_{h}))$, and hence $\sigma^{2d}=1$ in
$\Gamma_2(B_n(\Sigma_{h}))/\Gamma_3(B_n(\Sigma_{h}))$ (recall that,  by \reth{gam3},
$\Gamma_2(B_n(\Sigma_{h}))/\Gamma_3(B_n(\Sigma_{h}))=
\ang{\sigma^2}\cong \Z_{h+n-1}$). But this would imply
that $h+n-1\leq d$~--~a contradiction. This proves
part~(\ref{it:g2sigma}).

Part~(\ref{it:g3sigma}) may be proved in the same way
as~(\ref{it:gam3g}) of \reth{gam3}; indeed, the quotient $B_n(\Sigma_{g,m})/
\Gamma_3(B_n(\Sigma_{g,m}))$ has a presentation similar to that
of~(\ref{eq:presbngam4}) 
(with the $z_i$ central, but without the last relation), and is
isomorphic to $B_n(\Sigma_{g,m})/H$, 
where $H$ is the normal closure of $\sigma_1\sigma_2^{-1}$ in
$B_n(\Sigma_{g,m})$, and thus is equal 
to the normal closure of $\Gamma_2(B_n)$ in $B_n(\Sigma_{g,m})$. As in
\reth{gam3}, one may show 
that $H=\Gamma_3(B_n(\Sigma_{g,m})$ is perfect for $n\geq 5$.

Finally, to prove part~(\ref{it:resid}), as in \reth{gam3} it
suffices to prove that $\Gamma_3(B_n(\Sigma_{g,m}))\neq \brak{1}$. 
Suppose that $\Gamma_3(B_n(\Sigma_{g,m}))=\brak{1}$. Then
$\Gamma_2(B_n(\Sigma_{g,m}))\cong\Z$ by~(\ref{it:g2sigma}), 
and since $B_n(\Sigma_{g,m})\supset B_n$, it follows that $\Gamma_2(B_n)$ is cyclic; but since $n\geq 3$,
this contradicts the results of~\cite{GL}.
\end{prooftext}

%%%%%%%%%%%%%%%%%%%%%%%%%%%%%%%%%%%%%%%%%%%%%%%%%%%%%%%%%%%%%%%%%%%%%%%
%%%%%%%%%%%%%%%%%%%%%%%%%%%%%%%%%%%%%%%%%%%%%%%%%%%%%%%%%%%%%%%%%%%%%%%
%%%%%%%%%%%%%%%%%%%%%%%%%%%%%%%%%%%%%%%%%%%%%%%%%%%%%%%%%%%%%%%%%%%%%%%
%%%%%%%%%%%%%%%%%%%%%%%%%%%%%%%%%%%%%%%%%%%%%%%%%%%%%%%%%%%%%%%%%%%%%%%
%%%%%%%%%%%%%%%%%%%%%%%%%%%%%%%%%%%%%%%%%%%%%%%%%%%%%%%%%%%%%%%%%%%%%%%
%%%%%%%%%%%%%%%%%%%%%%%%%%%%%%%%%%%%%%%%%%%%%%%%%%%%%%%%%%%%%%%%%%%%%%%
%%%%%%%%%%%%%%%%%%%%%%%%%%%%%%%%%%%%%%%%%%%%%%%%%%%%%%%%%%%%%%%%%%%%%%%
%%%%%%%%%%%%%%%%%%%%%%%%%%%%%%%%%%%%%%%%%%%%%%%%%%%%%%%%%%%%%%%%%%%%%%%
%%%%%%%%%%%%%%%%%%%%%%%%%%%%%%%%%%%%%%%%%%%%%%%%%%%%%%%%%%%%%%%%%%%%%%%

\section{Braid groups on $2$ strands: properties and open questions}\label{sec:b2t2}

The aim of this section is to prove \reth{renil}. Consider first the
group presentation given by \reth{presbng}, and take $n=2$ and
$g=1$.   Setting $\alpha=a \sigma_1$, $\beta= b \sigma_1$ and $\gamma=
a\sigma_1 b$, one obtains the following presentation of $B_2(\T)$:
\begin{thm}[\cite{BG}]\label{th:alpha}
$B_2(\T)$ is generated by $\alpha$, $\beta$ and $\gamma$, subject to
the relations:
\begin{gather*}
\text{$\alpha^2$ and $\beta^2$ are central}\\
\alpha^2\beta^2=\gamma^2.
\end{gather*}
Further, $\alpha^2$ and $\beta^2$ generate the centre of $B_2(\T)$.
\end{thm}

Let $\map{p}{B_2(\T)}[B_2(\T)/Z(B_2(\T))]$ denote the canonical
projection. From this presentation, it follows that
$B_2(\T)/Z(B_2(\T))$ is generated by $\overline{\alpha}=p(\alpha)$,
$\overline{\beta}=p(\beta)$ and $\overline{\gamma}=p(\gamma)$, subject
to the relations $\overline{\alpha}^2= \overline{\beta}^2=
\overline{\gamma}^2=1$. So $B_2(\T)/Z(B_2(\T))$, which we identify
with $\Z_2 \ast \Z_2 \ast \Z_2$, is the Coxeter group
$W(\overline{\alpha}, \overline{\beta}, \overline{\gamma})$ associated
to the free group $\F[3](\overline{\alpha}, \overline{\beta},
\overline{\gamma})$, and $B_2(\T)$ is a central extension of
$W(\overline{\alpha}, \overline{\beta}, \overline{\gamma})$:
\begin{equation*}
1\to Z(B_2(\T)) \to B_2(\T)  \stackrel{p}{\to} \Z_2 \ast \Z_2 \ast
\Z_2\to 1.
\end{equation*}

This presentation of $B_2(\T)$ was considered in~\cite{BG}, where the
following length functions $\ell_{\widehat{\alpha}},
\ell_{\widehat{\beta}}$ were defined.  If $x\in
\brak{\alpha,\beta,\gamma}$, set:
\begin{equation*}
\ell_{\widehat{\alpha}}(x)=
\begin{cases}
& \text{$1$ if $x\neq \alpha$}\\
& \text{$0$ if $x=\alpha$,}
\end{cases}
\end{equation*}
and similarly for $\ell_{\widehat{\beta}}$. From \reth{alpha}, it
follows that each 
of $\ell_{\widehat{\alpha}}$ and $\ell_{\widehat{\beta}}$ extends to a homomorphism of $B_2(\T)$ onto $\Z$.

The following observation will be used in the proof of \reth{renil}.

\begin{prop}\label{center}
The intersection of $\Gamma_2(B_2(\T))$ and $Z(B_2(\T))$ is trivial.
\end{prop}

\begin{proof}
Let $x\in Z(B_2(\T))$. By \reth{alpha}, there exist $m,n\in\Z$ such
that $x=a^{2m}b^{2n}$, 
and thus $\ell_{\widehat{\alpha}}(x)=2n$ and
$\ell_{\widehat{\beta}}(x)=2m$. 
But $x\in\Gamma_2(B_2(\T))$, so $\ell_{\widehat{\alpha}}(x)=
\ell_{\widehat{\beta}}(x)=0$.
 We conclude that $m=n=0$, and hence $x=1$.
\end{proof}

We are now able to prove \reth{renil}.

\begin{prooftext}{\reth{renil}}\mbox{}
Set $G=\Z_2 \ast \Z_2 \ast \Z_2$.
\begin{enumerate}[(a)]
\item Suppose that $x\in \bigcap_{i\in\N}\; \Gamma_i(B_2(\T))$. Then
$p(x)\in \bigcap_{i\in\N}\; \Gamma_i(G)$, but
since $G$  is residually nilpotent~\cite{G}, it follows that $x\in \ker p=Z(B_2(\T))$. So $x=1$ by
Proposition~\ref{center}, and hence $B_2(\T)$ is residually nilpotent.
\item 
\begin{enumerate}[(i)]
\item Let us consider the following commutative diagram of short exact sequences:
\begin{equation*}
\xymatrix{%
1 \ar[r] & \Gamma_2(B_2(\T)) \ar[r] \ar[d]_{p_2} &
B_2(\T) \ar[d]_{p}  \ar[r] & B_n(\T)/\Gamma_2(B_n(\T)) \ar[d] \ar[r] & 1\\
1 \ar[r] & \Gamma_2(G) \ar[r] &
G  \ar[r] & G/\Gamma_2(G) \ar[r] & 1
}
\end{equation*}
The first and third vertical arrows are those induced by $p$. More
generally, for $i\geq 2$, let
$\map{p_i}{\Gamma_i(B_2(\T))}[\Gamma_i(G)]$ denote the epimorphism
induced by $p$. But it follows from Proposition~\ref{center} that
$p_2$ is also injective, so is an isomorphism. Since for $i\geq 3$,
$p_i$ is the restriction of $p_2$ to $\Gamma_i(B_2(\T))$, $p_i$ is an
isomorphism too. 

\item From~(\ref{it:igeq2})(\ref{it:gammai}), it follows that
 $\Gamma_i(B_2(\T))/\Gamma_{i+1}(B_2(\T)) \cong
 \Gamma_i(G)/\Gamma_{i+1}(G)$,
so it suffices to prove the result for $G$. We break the proof down into two parts as follows:
\begin{enumerate}[(1).]
\item Recall that the elements $\overline{\alpha}, \overline{\beta}$
  and $\overline{\gamma}$ are each of 
order~$2$, and generate $G$. We claim that every non-trivial element
  of $\Gamma_i(G)/\Gamma_{i+1}(G)$ is of 
order~$2$.  Since $\Gamma_i(G)/\Gamma_{i+1}(G)$ is a
  finitely-generated Abelian group by~\cite{MKS}, 
this will imply that it is isomorphic to a finite number, $R_i$ say,
  of copies of $\Z_2$. 
To prove the claim, recall from~\cite{MKS} that
  $\Gamma_i(G)/\Gamma_{i+1}(G)$ is  generated by 
the cosets modulo $\Gamma_{i+1}(G)$ of the $i$-fold simple commutators
  $[[\cdots[[\rho_1,\rho_2],\rho_3] \cdots , \rho_{i-1}], \rho_i]$, 
where $\rho_j\in\brak{\overline{\alpha}, \overline{\beta},
  \overline{\gamma}}$ for all $1\leq j\leq i$. We argue by induction
  on $i\geq 2$. 
Firstly, let $i=2$. Then $\Gamma_2(G)/\Gamma_3(G)$ is generated by the
  cosets of the $[\rho_1, \rho_2]$. 
But modulo $\Gamma_3(G)$, $[\rho_1,\rho_2]^2\equiv [\rho_1^2,
  \rho_2]\equiv 1$, 
and since $\Gamma_2(G)/\Gamma_3(G)$ is Abelian, this implies that all
  of its non-trivial elements are of order~$2$. 
Now suppose that $i\geq 3$, and suppose by induction that the result
  holds for $i-1$, so that $x^2\equiv 1$ modulo $\Gamma_i(G)$
 for all $x\in \Gamma_{i-1}(G)$. Every $i$-fold simple commutator may
  be written in the form $[x,\rho_i]$, where $x$ is a $(i-1)$-fold
  simple commutator, 
so belongs to $\Gamma_{i-1}(G)$, and
  $\rho_i\in\brak{\overline{\alpha}, \overline{\beta},
  \overline{\gamma}}$, so belongs to $G$. 
By the induction hypothesis, $x^2\in \Gamma_i(G)$, so $[x,\rho_i]^2
  \equiv [x^2, \rho_i]\equiv 1$ modulo $\Gamma_{i+1}(G)$, 
and once more, since $\Gamma_i(G)/\Gamma_{i+1}(G)$ is Abelian, all of
  its non-trivial elements are of order~$2$. 
This proves the claim.

\item The number $R_i$ of summands of $\Z_2$ is given by Theorem~3.4
  of~\cite{Ga}. We refer to Gaglione's notation in what follows. 
Since $U_{\infty}(x)=0$, the $R_{\infty}^j$ are all zero
  ($R_{\infty}^j$ represents the rank of the free Abelian 
factor of $\Gamma_j(G)/\Gamma_{j+1}(G)$), and so $R_i$ is as given in
  the statement of the theorem. 
It just remains to determine $k\alpha_k$ for all $k\geq 2$. A simple
  calculation shows that $1-U(x)=(1+x)^2(1-2x)$, 
hence:
\begin{equation*}
\frac{d}{dx} \ln \left( 1-U(x) \right)= \frac{2}{x+1}+\frac{2}{2x-1},
\end{equation*}
and that for $k\geq 2$, 
\begin{equation*}
\frac{d^k}{dx^k} \ln \left( 1-U(x) \right)= (-1)^{k+1}(k-1)! \left(\frac{2}{(x+1)^k}+
\frac{2^k}{(2x-1)^k}\right).
\end{equation*}
So 
\begin{equation*}
k\alpha_k=-\frac{1}{(k-1)!} \left. \left(\frac{d^k}{dx^k} \ln \left( 1-U(x) \right)\right) \right\vert_{x=0}=
2^k+2(-1)^k,
\end{equation*}
as required.
\end{enumerate}
\end{enumerate}
\item Given a group $G$, the quotient group $D_i(G)/D_{i+1}(G)$ is
torsion free and it is isomorphic to  $\Gamma_i(G)/\Gamma_{i+1}(G)$
modulo torsion, for $i\ge 1$ \cite{P}. Therefore, from part~(b) one deduces that
$D_2(B_2(\T))=D_{3}(B_2(\T))$.   On the other hand one can easily verify that
$B_2(\T)/\Gamma_{2}(B_2(\T)) \cong \Z^2 \oplus \Z_2$ and therefore
$B_2(\T)/D_2(B_2(\T))\cong \Z^2$. Since  $B_2(\T)$ is not Abelian, it follows
that $D_2(B_2(\T))$ is not trivial and then $\bigcap_{i\in\N} D_i(B_2(\T)\neq\brak{1}$.
\end{enumerate}
\end{prooftext}

\begin{rems}
From \reth{renil} one concludes that
$\Gamma_i(B_2(\Sigma_{1,p}))\neq \Gamma_{i+1}(B_2(\Sigma_{1,p}))$. On
the other hand, the group $\Gamma_2(B_2(\Sigma_{g,p}))$ is generated
by the set of conjugates of the commutators of the form $[g,g']$,
where $g,g'$ are generators of $B_2(\Sigma_{g,p})$. Therefore
$\Gamma_2(B_2(\Sigma_{g,p})) \subset P_2(\Sigma_{g,p})$. Since 
$P_2(\Sigma_{g,p})$ is residually  nilpotent for $p\ge 1$ (see
Section~\ref{puresection}), one deduces that $B_2(\Sigma_{g,p})$ is
residually soluble for $p\ge 1$. The question of whether
$B_2(\Sigma)$ is in fact residually nilpotent, when
$\Sigma$ is a surface of positive genus possibly with boundary different from the torus,
is open.
\end{rems}

To finish this section, we prove the following result:
\begin{prop} \label{biord}
The group $B_2(\T)$ is not bi-orderable.
\end{prop}

\begin{proof}
Consider $B_2(\T)$ with the group presentation in \reth{alpha}.
Set $g=\alpha \beta \gamma^{-1}$. The following equality holds in $B_2(\T)$:
\begin{equation*}
((\alpha \gamma)^{-1}  g (\alpha \gamma)) (\gamma^{-1} g \gamma) (\alpha^{-1}  g \alpha) (g) =1\, .
\end{equation*}
Since $g \not=1$, the group $B_2(\T)$ is not bi-orderable by Proposition~\ref{test}.
\end{proof}

Let $\Sigma$ be an orientable surface, possibly with boundary. If $n\ge
3$, $B_n(\Sigma)$ is not bi-orderable since it contains a copy of
$B_n$ which is not bi-orderable~\cite{Go}. If $n=1$, the group
$B_1(\Sigma)$ is isomorphic to $\pi_1(\Sigma)$ which is known to be
residually free. Therefore it is also residually torsion-free
nilpotent and hence bi-orderable.

\begin{rem}
If $\Sigma$ is an orientable surface, possibly with boundary,
different from the torus, the sphere and  the disc, the question of whether $B_2(\Sigma)$ is in fact
bi-orderable is
open. 
\end{rem}

%%%%%%%%%%%%%%%%%%%%%%%%%%%%%%%%%%%%%%%%
%%%%%%%%%%%%%%%%%%%%%%%%%%%%%%%%%%%%%%%%
%%%%%%%%%%%%%%%%%%%%%%%%%%%%%%%%%%%%%%%%
%%%%%%%%%%%%%%%%%%%%%%%%%%%%%%%%%%%%%%%%
%%%%%%%%%%%%%%%%%%%%%%%%%%%%%%%%%%%%%%%%
%%%%%%%%%%%%%%%%%%%%%%%%%%%%%%%%%%%%%%%%
%%%%%%%%%%%%%%%%%%%%%%%%%%%%%%%%%%%%%%%%
%%%%%%%%%%%%%%%%%%%%%%%%%%%%%%%%%%%%%%%%
%%%%%%%%%%%%%%%%%%%%%%%%%%%%%%%%%%%%%%%%
%%%%%%%%%%%%%%%%%%%%%%%%%%%%%%%%%%%%%%%%

\section{Residual torsion free nilpotence of surface pure braid groups}\label{puresection}

In this section we give a short survey on relations between surface
braids and mapping classes, and we show that pure braid groups of
surfaces with non-empty boundary may be realised as subgroups of
Torelli groups of surfaces with one boundary component. 

%%%%%%%%%%%%%%%%%%%%%%%%%%%%%%%%%%%%%%%
%%%%%%%%%%%%%%%%%%%%%%%%%%%%%%%%%%%%%%%
%%%%%%%%%%%%%%%%%%%%%%%%%%%%%%%%%%%%%%%
%%%%%%%%%%%%%%%%%%%%%%%%%%%%%%%%%%%%%%%

\subsection{Surface pure braid groups}

We start by recalling a group presentation for pure braid groups of surfaces with one boundary
component~\cite{B}.

\begin{thm} \label{th:purepres}
Let  $\Sigma_{g,1}$ be a compact, connected orientable surface 
of genus $g\ge 1$ with one boundary component. The group $P_n(\Sigma_{g,1})$ admits the following
presentation:
\begin{enumerate}
\item[\textbf{Generators:}] $\{A_{i,j}\; | \;1 \le i \le 2g+ n -1,  2g +1\le j \le 2g + n, i<j \}.$
\item[\textbf{Relations:}] 
\begin{eqnarray*}
 &\text{(PR1)}&  A_{i,j}^{-1}  A_{r,s} A_{i,j} = A_{r,s} \;  \; \mbox{if} \, \,(i<j<r<s)  \;   \mbox{or} \,
(r+1<i<j<s),\\ 
 &     &  \mbox{or} \, (i=r+1<j<s \, \, \mbox{for even} \, \, r<2g  \, \, \mbox{or}   \, \, r>2g  ) \,; \\ 
 &\text{(PR2)}&  A_{i,j}^{-1}  A_{j,s} A_{i,j} = A_{i,s}  A_{j,s} A_{i,s}^{-1} \;  \; \mbox{if} \, \,
(i<j<s)\,;\\
 &\text{(PR3)}&  A_{i,j}^{-1}  A_{i,s} A_{i,j} = A_{i,s} A_{j,s} A_{i,s}  A_{j,s}^{-1} A_{i,s}^{-1} \; \;
\mbox{if} \, \, (i<j<s)\,; \\
 &\text{(PR4)}&
A_{i,j}^{-1}A_{r,s}A_{i,j}=A_{i,s}A_{j,s}A_{i,s}^{-1}A_{j,s}^{-1}A_{r,s}A_{j,s}A_{i,s}A_{j,s}^{-1}A_{i,s}^{-1}
\\
 &     &   \mbox{if} \, \,(i+1<r<j<s)  \;   \mbox{or} \\
 &     &   \, \, (i+1=r<j<s  \; \mbox{for odd }  \, \, r<2g  \, \, \mbox{or}   \, \, r>2g ) \,; 
 \end{eqnarray*}
 \begin{eqnarray*}
 &\text{(ER1)}&  A_{r+1,j}^{-1} A_{r,s} A_{r+1,j}=A_{r,s} A_{r+1,s} A_{j,s}^{-1} A_{r+1,s}^{-1}  \\ 
 &     &     \mbox{if} \,  \,r \, \mbox{odd and}\, \,r<2g                 \,               ; \\
 &\text{(ER2)}&   A_{r-1,j}^{-1}  A_{r,s} A_{r-1,j}=  A_{r-1,s} A_{j,s} A_{r-1,s}^{-1}  A_{r,s} A_{j,s}
A_{r-1,s} A_{j,s}^{-1}A_{r-1,s}^{-1} \\ 
 &    &   \mbox{if} \,  \,r \, \mbox{even and}\, \,r<2g   \,               . \\
 \end{eqnarray*}
 \end{enumerate}
\end{thm}

\begin{figure}[h]
\begin{center}
\psfig{figure=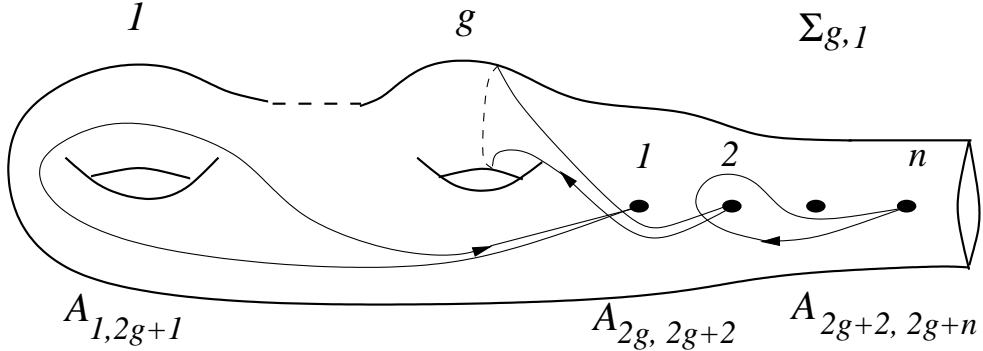,width=13cm}%,height=2cm}}
\caption{\label{generateur}Projection of representatives of the generators $A_{i,j}$. We
represent $A_{i,j}$ by its only non-trivial strand.}
\end{center}
\end{figure}

As a representative of the generator $A_{i,j}$, we may take a
geometric braid whose only non-trivial (non-vertical) strand is the $(j-2g)$th
one.  In Figure~\ref{generateur}, we illustrate the projection of such braids on the surface $\Sigma_{g,1}$
(see also Figure~8 of~\cite{B}). Some misprints in Relations~(ER1) and~(ER2) of Theorem~5.1 of~\cite{B} have
been corrected.

%%%%%%%%%%%%%%%%%%%%%%%%%%%%%%%%%%%%%%%%%%%%%%%%%%%%%%%%%%%%%%%%%%%%%%%%%%%%%%%%
  
With respect to the presentation of $B_n(\Sigma_g)$ given
in  \reth{presbng}, the elements $A_{i,j}$ are the following braids:
\begin{itemize}
\item $A_{i,j}= \sigma_{j-2g} \cdots  \sigma_{i+1-2g} \sigma_{i-2g}^2 \sigma_{i+1-2g}^{-1}
 \cdots \sigma_{j-2g}^{-1}$, for $i \ge 2g$;
\item $A_{2i,j}= \sigma_{j-2g} \cdots  \sigma_{1}  a_{g-i+1} \sigma_{1}^{-1} \cdots \sigma_{j-2g}^{-1}$, for 
$1 \le i \le g$;
\item $A_{2i-1,j}= \sigma_{j-2g} \cdots  \sigma_{1}  b_{g-i+1} \sigma_{1}^{-1} \cdots \sigma_{j-2g}^{-1}$,
for 
$1 \le i \le g$.
\end{itemize}

Relations~(PR1), \dots,~(PR4) correspond to the classical relations
for the pure braid group $P_n$~\cite{Bir}. New relations arise when
we consider two generators $A_{2i,j}$, $A_{2i-1,k}$, for $1 \le i \le
g$ and $j \not= k$. They correspond to loops based at two
different points which go around the same handle.

%%%%%%%%%%%%%%%%%%%%%%%%%%%%%%%%%%%%%%%
%%%%%%%%%%%%%%%%%%%%%%%%%%%%%%%%%%%%%%%
%%%%%%%%%%%%%%%%%%%%%%%%%%%%%%%%%%%%%%%
%%%%%%%%%%%%%%%%%%%%%%%%%%%%%%%%%%%%%%%

\subsection{Mapping class groups, bounding pair braids and pure braids}\label{sec:bound}

The \emph{mapping class group} of a surface $\Sigma_{g,p}$, denoted by $\mathcal{M}_{g,p}$, is the group of
isotopy classes of orientation-preserving self-homeomorphisms  which
fix the boundary components pointwise. If the surface has empty boundary then we shall just write
$\MG$. Note that we will denote the composition in
the mapping class groups from left to right\footnote{We do this in
order to have the same group-composition in braid groups and mapping
class groups.}.

Let $\mathcal{P}=\{x_1, \dots, x_n\}$ be a set of $n$ distinct points
in the interior of the surface $\Sigma_{g,p}$.  The \emph{punctured
mapping class group} of $\Sigma_{g}$ relative to $\mathcal{P}$ is
defined  to be the group of isotopy classes of orientation-preserving
self-homeomorphisms which fix the boundary components pointwise, and
which fix $\mathcal{P}$ setwise. This group, denoted by $\PMG$, does
not depend on the choice of  $\mathcal{P}$, but just on its cardinal.
We define the \emph{pure punctured mapping class group}, denoted by $\PPMG$, to be the
subgroup of (isotopy classes of) homeomorphisms which fix the set
$\mathcal{P}$ pointwise. 

We note $T_C$ a Dehn twist along a simple closed curve $C$. 
Let $C$ and $D$ be two simple closed curves bounding an annulus
containing the single puncture $x_j$.  We shall say that the
multitwist $T_{C} T_{D}^{-1}$ is a \emph{$j$-bounding pair braid}.

%%%%%%%%%%%%%%%%%%%%%%%%%%%%%%%%%%%%%%%
%%%%%%%%%%%%%%%%%%%%%%%%%%%%%%%%%%%%%%%
%%%%%%%%%%%%%%%%%%%%%%%%%%%%%%%%%%%%%%%
%%%%%%%%%%%%%%%%%%%%%%%%%%%%%%%%%%%%%%%

%\subsection{Pure braids, bounding pair braids and Torelli groups}\label{sec:bound}

Surface braid groups are related to mapping class groups as follows:

\begin{thm}{\bf (Birman \cite{Bir})} \label{th:bir}
Let  $g\ge 1$ and  $p \ge 0$.  Let $\psi: \PMG \to \mathcal{M}_{g,p}$
be the homomorphism induced by the map which forgets the set
$\mathcal{P}$. If $\Sigma_{g,p}$ is  different from  the torus  then 
$\ker \psi$ is isomorphic to  $B_n(\Sigma_{g,p})$.
\end{thm}

\begin{rem}\label{rem:pure}
In particular, if  $\Sigma$ is an orientable surface (possibly with
boundary) of positive genus and different from the torus, the surface
pure braid group $P_n(\Sigma)$ may  be identified with the subgroup of $\PMG$
generated by bounding pair braids (see for instance~\cite{Bir}, where bounding pair braids are called
spin-maps).  
\end{rem}

%%%%%%%%%%%%%%%%%%%%%%%%%%%%%%%%%%%%%%%
%%%%%%%%%%%%%%%%%%%%%%%%%%%%%%%%%%%%%%%
%%%%%%%%%%%%%%%%%%%%%%%%%%%%%%%%%%%%%%%
%%%%%%%%%%%%%%%%%%%%%%%%%%%%%%%%%%%%%%%

\subsection{Torelli groups}

We recall that the Torelli group $\mathcal{T}_{g,1}$ is the subgroup
of the mapping class group $\mathcal{M}_{g,1}$ which acts trivially on
the first homology group of the surface $\Sigma_{g,1}$.

Before stating the main theorem of this section, we recall the following exact sequence:
\begin{equation}\label{eq:sequence}
1 \to \Z^n  \to \mathcal{M}_{g,n+p}
  \stackrel{q}{\to} \PPMG \to 1 \,,
\end{equation}
where $\Z^n$ is central and generated by Dehn twists along the first $n$ boundary
components of  $\Sigma_{g,n+p}$. Geometrically, the projection $q$ may
be obtained by glueing one-punctured discs $\D_{1},\ldots,\D_{n}$,
say, onto the first $n$  boundary components. 
Note that sequence~\ref{eq:sequence} does not split, since the first homology group of $\mathcal{M}_{g,n+p}$ is 
trivial when $g$ is greater than 2. Nevertheless, as explained in Lemma~\ref{lem:puretorelli}, when the
surface has 
boundary, we have a section over the corresponding pure braid group.

\begin{lem}\label{lem:puretorelli}
Let $\Sigma_{g,1}$  be a surface of genus greater than or equal to  one with one
boundary component. Then the group $P_n(\Sigma_{g,1})$ embeds
in  $\mathcal{T}_{g+n,1}$.
\end{lem}

\begin{proof}
Applying \reth{bir} and Remark~\ref{rem:pure}, we identify
$P_n(\Sigma_{g,1})$ with the subgroup of $\PPMGa$  generated by
bounding pair braids. Let us first embed $P_n(\Sigma_{g,1})$ in
$\mathcal{M}_{g,n+1}$. To achieve this, we construct a section $s$ on
$P_n(\Sigma_{g,1})$ of the sequence~(\ref{eq:sequence}). For each
generator $A_{i,j}$ of $P_n(\Sigma_{g,1})$, we define  $s(A_{i,j})$ as
follows. Consider two simple closed  curves $a$ and $a'$ lying in
$\Sigma_{g,1}$ such that $A_{i,j}$ is equal to the boundary pair
braids  $T_{a}T_{a'}^{-1}$. These two curves may be chosen so as to
avoid the discs  $\D_{1},\ldots,\D_{n}$, and thus may be seen as lying
in $\Sigma_{g,n+1}$. If $d_{j}$ is a simple closed curve  parallel to
the $j$th-boundary component, we set 
$s(A_{i,j})=T_{a}T_{a'}^{-1}T_{d_{j}}$, which we denote by $A'_{i,j}$.
Since the Dehn twists $T_{d_{1}},\ldots, T_{d_{n}}$ belong to the
kernel of  $q$, one has $q\circ s=\mathrm{Id}$, and hence $s$ is
injective. We claim that $s$ is a homomorphism. To prove this, we have
to show that relations~(PR1-4) and~(ER1-2) are satisfied in
$\mathcal{M}_{g,n+1}$ via\nolinebreak[4] $s$.

The four first relations may be written in the form
$hA_{r,s}h^{-1}=A_{r,s}$, where $h$ is a word in the  $A_{i,j}$'s.
These relations are compatible with $s$, since for all simple closed
curves $a$ in $\Sigma_{g,n+1}$, and all $h$ in  $\mathcal{M}_{g,n+1}$,
one has: 
\begin{equation}\label{eq:relation}
T_{h(a)}=h^{-1}T_{a}h.
\end{equation}

\begin{figure}
\begin{center}
%\scalebox{0.5}
{\psfig{figure=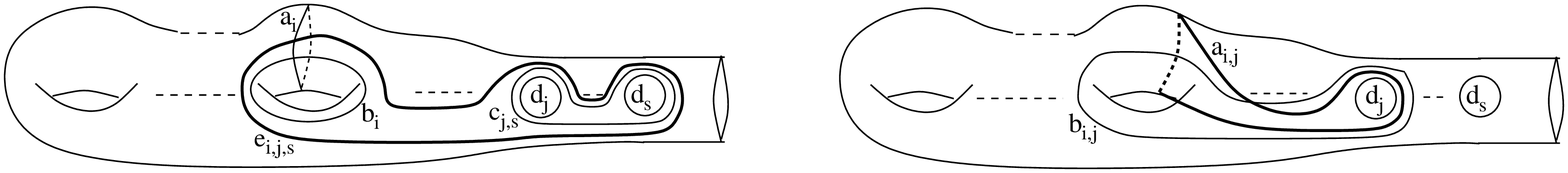,width=14cm,height=2cm}}
\caption{\label{courbes}curves for relations~(PR2) and~(ER2)}
\end{center}
\end{figure}  

For example, relation~(PR1) is compatible with $s$ because the curves occurring
in $A'_{i,j}$ are disjoint from those occuring in $A'_{r,s}$.
For~(PR2), the bounding pair braid $A_{j,s}$ (resp.\ $A_{i,j}$,
$A_{i,s}$) is equal to  $T_{d_{j}}T_{c_{j,s}}^{-1}$ (resp.\
$T_{a}T_{a'}^{-1}$, $T_{b}T_{b'}^{-1}$ for 
$(a,a')\in\{(a_{i},a_{i,j}),(b_{i},b_{i,j}), (d_{i},c_{i,j})\}$
and  $(b,b')\in\{(a_{i},a_{i,s}), (b_{i},b_{i,s}),(d_{i},c_{i,s})\}$,
where  curves are those described by Figure~\ref{courbes}).  Thus
we have:
\[ \begin{array}{rcl}
A'^{-1}_{i,j}A'_{j,s}A'_{i,j} & = &
[T^{-1}_{a}T_{a'}T^{-1}_{d_{j}}]T_{d_{j}}T^{-1}_{c_{j,s}}T_{d_{s}}[T_{d_{j}}T_{a'}^{-1}T_{a}] \\
 & = & [T^{-1}_{a}T_{a'}]T^{-1}_{c_{j,s}}[T_{a'}^{-1}T_{a}]T_{d_{j}}T_{d_{s}} \ \ \hbox{ since the }\ 
 T_{d_{k}}\mathrm{'s}\ 
 \hbox{ are central}\\
 & = & T^{-1}_{T_{a}T^{-1}_{a'}(c_{j,s})}T_{d_{j}}T_{d_{s}}\ \ \hbox{ by~(\ref{eq:relation})},
\end{array}
\]
and similarly, $A'_{i,s}A'_{j,s}A'^{-1}_{i,s}=T^{-1}_{T^{-1}_{b}T_{b'}(c_{j,s})}T_{d_{j}}T_{d_{s}}$. Now, it 
is easy to see that
$$T_{a}T^{-1}_{a'}(c_{j,s})=T^{-1}_{a'}(c_{j,s})=T^{-1}_{b}T_{b'}(c_{j,s}),$$
which yields the required relation. The compatibility of relations~(PR3-4) with $s$ may be proved in the same
way; we leave this as an exercise for the reader.

Relation~(ER1) is a consequence of the lantern relation and  relation~(\ref{eq:relation}). Indeed, if we
consider the seven curves $b_{i}$, $d_{j}$, $d_{s}$, $e_{i,j,s}$, $b_{i,s}$, 
$b_{j,s}$ and $c_{j,s}$ shown in Figure~\ref{courbes} (where $r=2i-1$), the lantern relation 
may be written as:
$$T_{e_{i,j,s}}T_{b_{i}}T_{d_{j}}T_{d_{s}}=T_{b_{i,s}}T_{b_{i,j}}T_{c_{j,s}},$$
which implies that
$$A'_{r,s}=T_{b_{i}}T^{-1}_{b_{i,s}}T_{d_{s}}=T_{b_{i,j}}T^{-1}_{e_{i,j,s}}T_{c_{j,s}}T_{d_{j}}^{-1}=
T_{b_{i,j}}T^{-1}_{e_{i,j,s}}T_{d_{s}}A'^{-1}_{j,s}.$$
Since $A'_{r+1,j}=T_{a_{i}}T^{-1}_{a_{i,j}}T_{d_{j}}$, we obtain
\[\renewcommand{\arraystretch}{1.8}
\begin{array}{rcl}
A'^{-1}_{r+1,j}A'_{r,s}A'_{r+1,j} & = & 
\Bigl[T^{-1}_{d_{j}}T_{a_{i,j}}T^{-1}_{a_{i}}T_{b_{i,j}}T^{-1}_{e_{i,j,s}}T_{d_{s}}T_{a_{i}}T^{-1}_{a_{i,j}}T_{d_{j}}\Bigr]
\Bigl[A'^{-1}_{r+1,j}A'^{-1}_{j,s}A'_{r+1,j}\Bigr] \\
 & = & 
\Bigl[T_{a_{i,j}}T^{-1}_{a_{i}}T_{b_{i,j}}T^{-1}_{e_{i,j,s}}T_{a_{i}}T^{-1}_{a_{i,j}}\Bigr]T_{d_{s}}
\Bigl[A'_{r+1,s}A'^{-1}_{j,s}A'^{-1}_{r+1,s}\Bigl]\ \hbox{ by (PR2)}\\
& = & T_{T_{a_{i}}T^{-1}_{a_{i,j}}(b_{i,j})}T^{-1}_{T_{a_{i}}T^{-1}_{a_{i,j}}(e_{i,j,s})}T_{d_{s}}
\Bigl[A'_{r+1,s}A'^{-1}_{j,s}A'^{-1}_{r+1,s}\Bigl]\ \hbox{ by (\ref{eq:relation}). }
\end{array}
\renewcommand{\arraystretch}{1}\]
But $T_{a_{i}}T^{-1}_{a_{i,j}}(b_{i,j})=b_{i}$ and $T_{a_{i}}T^{-1}_{a_{i,j}}(e_{i,j,s})=b_{i,s}$, so
\[\renewcommand{\arraystretch}{1.8}
\begin{array}{rcl}
A'^{-1}_{r+1,j}A'_{r,s}A'_{r+1,j} & = & T_{b_{i}}T^{-1}_{b_{i,s}}T_{d_{s}}
\Bigl[A'_{r+1,s}A'^{-1}_{j,s}A'^{-1}_{r+1,s}\Bigl] \\
& = & A'_{r,s}A'_{r+1,s}A'^{-1}_{j,s}A'^{-1}_{r+1,s}\,,
\end{array}
\renewcommand{\arraystretch}{1}\]
which is~relation (ER1). Relation~(ER2) is also a consequence of a lantern: again, we leave the details to 
the reader.

Hence $s:P_n(\Sigma_{g,1})\to \mathcal{M}_{g,n+1}$ is an embedding. Glueing a one-holed torus onto each of
the first
$n$ boundary components of  $\Sigma_{g,n+1}$, we obtain a homomorphism \linebreak[4]$\phi:\mathcal{M}_{g,n+1}
\to 
\nolinebreak[4]\mathcal{M}_{g+n,1}$ which is injective (see~\cite{PR2}). Clearly, the image under $\varphi$
of each 
$s(A_{i,j})$ acts trivially on the homology group $H_{1}(\Sigma_{g+n,1};\Z)$. Thus $\varphi \circ
s(P_n(\Sigma_{g,1}))$ lies in the Torelli group of $\Sigma_{g+n,1}$.
\end{proof}

\begin{rem}\label{rem:closed}
This embedding of $P_n(\Sigma_{g,1})$ in $\mathcal{T}_{g+n,1}$ does not hold for surfaces with empty
boundary. Indeed, the 
group $P_n(\Sigma_{g})$ has an extra relation~(TR) (see Theorem 5.2 of~\cite{B}) which is not satisfied by
the section $s$; 
if $L$ (resp. $R$) is the left-hand side (resp. right-hand side) of this relation, one can check using 
lantern relations that we have $s(L)=s(R)d_{k}^{2(g-1)}$ in $\mathcal{M}_{g,n}$ ($k$ is the same index as in 
the relation~(TR) in Theorem 5.2 of \cite{B}).
Nevertheless, it would be interesting to know whether the sequence~(\ref{eq:sequence}) splits over the pure braid 
group $P_n(\Sigma_{g})$.
\end{rem}

\begin{thm}
Let $\Sigma$ be  the torus, or a surface of positive genus 
with non-empty boundary. Then the group $P_n(\Sigma)$ is residually torsion-free nilpotent.
\end{thm} 

\begin{proof}
Let $\Sigma_{g,1}$ be a surface of genus greater than or equal to
one, and with one boundary component. First, we remark that
\relem{puretorelli} and the residually torsion-free nilpotence of Torelli groups (see for instance section 14 
of~\cite{H}) imply that  $P_n(\Sigma_{g,1})$  is residually 
torsion-free nilpotent. Now let  $\Sigma_{g, p}$ be a surface with $p
> 1$ boundary components. The group $P_n(\Sigma_{g,p})$ may be
realised as the  subgroup of $P_{n+p-1}(\Sigma_{g,1})$  formed by the
braids whose first $p-1$ strands are vertical. Therefore
$P_n(\Sigma_{g,p})$ is  residually torsion-free  nilpotent.

The remaining case is that of the pure braid group on $n$ strands of 
$\T$.  From  \relem{toreholetore}, one deduces easily that
$D_i(P_n(\T))=D_i(P_{n-1}(\Sigma_{1,1}))$ for $i>1$, and thus 
the group $P_n(\T)$ is residually torsion-free nilpotent.
\end{proof}

\begin{lem}\label{lem:toreholetore}
The group $P_n(\T)$ is isomorphic to $P_{n-1}(\Sigma_{1,1}) \times \Z^2$.
\end{lem}

\begin{proof}
Consider the pure braid group exact sequence for an orientable surface $\Sigma$:
$$
\quad 1 \to P_{n-1}(\Sigma \setminus \{x_1\}) \to
P_n(\Sigma) \stackrel{\theta}{\to} \pi_{1}(\Sigma) \to 1,
$$
where geometrically, $\theta$ is the map that forgets the paths
pointed at $x_2, \ldots, x_n$.  Since  $ZP_n(\Sigma_{1,1})$ is
trivial~\cite{PR1}, we deduce that the restriction of $\theta$ to
$ZP_n(\T)$ is injective. Since $ZP_n(\T)=\Z^2$~\cite{PR1}, the
restriction of $\theta$ to $ZP_n(\T)$ is in fact  an isomorphism, and we conclude that   $P_n(\T)$ is
isomorphic to the direct
product  $P_{n-1}(\Sigma_{1,1}) \times \Z^2$.
\end{proof}

\begin{rem}
In the case of the sphere, the group $P_n(\mathbb{S}^2)$ is
isomorphic to $\Z_2 \times P_{n-2}(\Sigma_{0,3})$ (see \cite{GG1}).
Therefore, for $i>1$, $\Gamma_i(P_n(S^2))$ and
$\Gamma_i(P_{n-2}(\Sigma_{0,3}))$ are isomorphic. Since  
$P_{n-2}(\Sigma_{0,3})$ is a subgroup of $P_{n}$ (which may be
realised geometrically as the subgroup of braids whose last strand is
vertical), from~\cite{FR} it follows that $P_n(\mathbb{S}^2)$ is  residually
nilpotent, but it is not residually torsion-free  nilpotent since
$P_n(\mathbb{S}^2)$ has torsion elements.
\end{rem}

\begin{rem}
Pure braid groups of surfaces of genus $g\ge 2$ with empty boundary are bi-orderable (\cite{Go}).  To the
best of our knowledge, it is not
known whether they are residually torsion-free nilpotent.
\end{rem}

%%%%%%%%%%%%%%%%%%%%%%%%%%%%%%%%%%%%%%%%
%%%%%%%%%%%%%%%%%%%%%%%%%%%%%%%%%%%%%%%%
%%%%%%%%%%%%%%%%%%%%%%%%%%%%%%%%%%%%%%%%
%%%%%%%%%%%%%%%%%%%%%%%%%%%%%%%%%%%%%%%%
%%%%%%%%%%%%%%%%%%%%%%%%%%%%%%%%%%%%%%%%
%%%%%%%%%%%%%%%%%%%%%%%%%%%%%%%%%%%%%%%%
%%%%%%%%%%%%%%%%%%%%%%%%%%%%%%%%%%%%%%%%
%%%%%%%%%%%%%%%%%%%%%%%%%%%%%%%%%%%%%%%%
%%%%%%%%%%%%%%%%%%%%%%%%%%%%%%%%%%%%%%%%
%%%%%%%%%%%%%%%%%%%%%%%%%%%%%%%%%%%%%%%%
%%%%%%%%%%%%%%%%%%%%%%%%%%%%%%%%%%%%%%%%
%%%%%%%%%%%%%%%%%%%%%%%%%%%%%%%%%%%%%%%%
%%%%%%%%%%%%%%%%%%%%%%%%%%%%%%%%%%%%%%%%

\end{document}